\documentclass{amsart}

\usepackage{graphicx}

\usepackage{caption}
\usepackage{subcaption}

\usepackage[pagewise]{lineno}

\newtheorem{theorem}{Theorem}[section]

\theoremstyle{definition}
\newtheorem{definition}[theorem]{Definition}

\theoremstyle{remark}

\numberwithin{equation}{section}

\begin{document}

\title{Autoencoding for the 'Good Dictionary' \\of eigen pairs of the Koopman Operator}

\author{Neranjaka Jayarathne}
\address{Center for Complex Systems Science, Clarkson University, Potsdam, NY 13699-5815}
\email{neranjaka.jayarathne@ieee.org}
\thanks{Neranjaka Jayarathne is supported from the ONR}

\author{Erik M. Bollt}
\address{Department of Electrical and Computer Engineering and The Clarkson Center for Complex Systems Science, Clarkson University, Potsdam, New York 13699, USA}
\email{ebollt@clarkson.edu}
\thanks{Erik M. Bollt is supported by the ONR, ARO and AFSOR and the NIH-CRCN}


\date{June 06, 2023}

\dedicatory{}


\begin{abstract}
Reduced order modelling relies on representing complex dynamical systems using simplified modes, which can be achieved through Koopman operator analysis. However, computing Koopman eigen pairs for high-dimensional observable data can be inefficient. This paper proposes using deep autoencoders, a type of deep learning technique, to perform non-linear geometric transformations on raw data before computing Koopman eigen vectors. The encoded data produced by the deep autoencoder is diffeomorphic to a manifold of the dynamical system, and has a significantly lower dimension than the raw data. To handle high-dimensional time series data, Takens's time delay embedding is presented as a pre-processing technique. The paper concludes by presenting examples of these techniques in action.
\end{abstract}

\maketitle

\section{Introduction}\label{ref_Introduction}
With the rapid advancement of computing machinery scientists acquired capability to analyse large data sets and extract information. Data driven techniques gathered momentum with this recent growth. The modern discipline of data driven sciences can be viewed as a combined outcome of available advanced mathematical techniques and high performance computing. According to \cite{brunton2022data} data driven discovery is revolutionizing the modeling, prediction and control in a diverse range of complex systems found in turbulence, the brain, climate, epidemiology, finance, robotics, and autonomy.

As opposed to analysing the behaviors of complex dynamical systems in terms of examining single trajectories, it has become possible to empirically discuss global questions in terms of evolution of density \cite{bollt2013applied}. Apart from the traditional analysis using geometrical techniques we can analyse the evolution of dynamical systems using the \textit{Frobenius-Perron} transfer operator \cite{bollt2013applied}. The left-adjoint of the \textit{Frobenius-Perror} operator is the \textit{Koopman} operator. \textit{Koopman} operator produces the values produced by a measurement function applied to the dynamical system. Instead of producing the value of the measurement it produces a value in future.

The data oriented \textit{Koopman} operator analysis perspective for analyzing dynamical systems has become extensively popular and relevant lately in science and engineering \cite{bollt2021geometric, budivsic2012applied,kutz2016dynamic,lan2013linearization}. This can be used to interpret dynamical systems irrespective of them being linear or non-linear \cite{bollt2021geometric}. However, the compromise is that the finite dimensional non-linear system may or may not require infinitely many dimensions to be represented as a linear system. For most of the computations requirements a small error could be tolerated. Therefore the infinite dimensionality can be truncated at the cost of reduced accuracy. In their work \cite{bollt2021geometric} describes the foundations for this work. Furthermore, it explains why the low dimensional dataset could be used to do the same. This paper brings out a machine learning technique to reduce the dimensionality of the dataset produced by a dynamical system using deep autoencoders. 

A \textit{good dictionary} in this context refers to a set of eigenvectors for an efficient representations of the dynamical system. When representing arbitrary observable functions by a series of eigenvectors, thee is freedom to define an efficient representation. Choosing the best from the rich and infinite eigenvectors develops a good dictionary of eigenvectors.  The reader is referred to \cite{avila2020data, mezic2005spectral, mezic2020spectrum, mezic2004comparison} for \textit{Koompan mood decomposition} as a global analysis of dynamical systems. This can be used to forcast as well as descriptive decompositions such as growing, fluctuating and decaying components. 

The utility of the initial DMD and exact DMD methods as numerical techniques for computing Koopman eigenpairs persists, as evidenced by their continued use in research \cite{bollt2021geometric, kutz2016dynamic, rowley2009spectral, schmid2010dynamic}. Several variants of DMD exist, including sparse DMD \cite{jovanovic2012low}, extended DMD (EDMD), and those that employ basis vectors to interpret data flow \cite{kevrekidis2016kernel, williams2015data}. In machine learning nomenclature, these basis vectors are commonly referred to as a \textit{dictionary of features} \cite{bollt2021geometric}. Following this dircetion \cite{li2017extended} describes an extended dynamic mode decomposition with optimal dictionary learning (EDMD-DL). DMD with control laws are introduced in \cite{kaiser2020data, kutz2016dynamic} which are emphasized on developing control laws.

In this paper we are developing a computationally efficient way of computing a good dictionary of Koopman eigen pairs. The technique of computing eigen pairs is based on \cite{bollt2021geometric}. Computational efficiency gain was induced by making a geometric transformation to the data manifold by means of deep learning techniques. In section \ref{ref_Mathematical Foundations} we present the mathematical and numerical foundations the Koopman eigen pairs are computed. In section \ref{ref_Manifold Learning} the machine learning technique was introduced. In section \ref{ref_Examples} examples are presented using the technique introduced. In section \ref{ref_Conclusions and Further Work} the conclusions are derived with a direction for future developments.

\section{Mathematical Foundations}\label{ref_Mathematical Foundations}

	\subsection{Koopman Analysis}
	First we review the underlying mathematics of Koopman analysis. Let,
	\begin{equation}
		\dot{X} = F(X), \quad F:M \to \mathbb{R}
	\end{equation} 
	be a dynamical system.
	The flow can be stated as a function for each $t \in \mathbb{R}$ such that $X(t) = \rho_t(X_0)$ where $X_0 = X(0) \in M$. This means that the trajectory starts at $X_0$ when $t=0$. A \textit{ Koopman Operator(KO))} describes the evolution of observables along the flow \cite{bollt2021geometric,budivsic2012applied, mezic2013analysis}. This study is based on analyzing these observables. Let us call these observables as observation functions. With reference to \cite{bollt2021geometric}, in this work we consider the observation function to be,
	\begin{equation}
		g: M \to \mathbb{C}
	\end{equation} where,
	\begin{equation}
		g \in L^2(M) = \left\{g: \int_{M}|g(s)|^2ds < \infty\right\}
	\end{equation}
	The scalar observation functions introduced here can be stacked to form vector observations. \textit{KO} defines how the observation functions evolve over time along the orbits of the dynamical system.
	
	\begin{definition}
		Let $\rho_t$ be a semiflow. That is, $\rho_t$ is define for $t \geq 0$. Then the \textit{KO} $\mathbb{K}_t$ is defined for $L^2(M)$ such that,
		\begin{equation}
			\mathbb{K}_t[g](X) = g \circ \rho_t(X)
		\end{equation}
		where, $\rho_t(X_0) = X_t$.
	\end{definition}

	In lehman's terms, for each $X$ we observe that value of an observable $g(x)$, not at $x$ itself, but after a time $t$ at $\rho_t(X)$. The $KO$ is linear on $L^2(M)$. This may be at the cost of it being infinite dimensional. \textit{eigenvectors} and \textit{eigenvalues} are studied under the spectral analysis of the $KO$ \cite{budivsic2012applied,gaspard2005chaos,mezic2013analysis}.
	
	\begin{definition}
		Suppose $\phi_{\lambda}$ is an eigenvector of $\mathbb{K}_t$ and $\lambda$ is the corresponding eigenvalue. Then the eigenpair should satisfy the equation,
		\begin{equation}
			\mathbb{K}_t[\phi_\lambda](X) = b^t\phi_\lambda(X)= e^{\lambda t}\phi_{\lambda}(X)
		\end{equation}
		or
		\begin{equation}
			\mathbb{K}_t [\phi_{\lambda}](X) = e^{\lambda t}\phi_{\lambda}(X).
		\end{equation}
	\end{definition}

	for each eigenvalue $\lambda$ there are uncountable many eigenvectors. In general these are not linearly independent. It is interesting to explore how an observable function can be decomposed into eigenvectors of the $KO$.
	
	\begin{definition}
		Let $\overline{g(X)}$ be a vector of observable functions $g_k(X)$.
		\begin{equation}
			\overline{g(X)} = \left[g_1(X), \dots, g_D(X)\right] : M \to \mathbb{C}^D.
		\end{equation}
	\end{definition}

	Suppose $\overline{g(X)}$ could be written as follows,
	
	\begin{equation}
		\overline{g(X)} = \sum_{j = 1}^{\infty} \phi_{\lambda_{j}}(X)\mathbf{v}_j.
	\end{equation}
	
	Then it follows that,
	\begin{equation}
		\overline{g(X_0)} = \sum_{j = 1}^{\infty} \phi_{\lambda_{j}}(X_0)\mathbf{v}_j.
	\end{equation}
	
	Apply $KO$ to oth sides of the equation,
	\begin{equation}
		\mathbb{K}\overline{g(X_0)} = \mathbb{K}\sum_{j = 1}^{\infty} \phi_{\lambda_{j}}(X_0)\mathbf{v}_j
	\end{equation}
	\begin{equation}
		g \circ \rho_t(X_0) = \sum_{j = 1}^{\infty} \mathbb{K}\left[\phi_{\lambda_{j}}\right](X_0)\mathbf{v}_j
	\end{equation}
	
	where $X_t = X(t)$.
	
	Since eigenvectors are not unique, the authors are interested in finding an efficient representation using eigenvectors. 
	\begin{definition}
		In \cite{bollt2021geometric} a \textbf{k}-efficient finite set of Koopman eigenpairs $\left\{\psi_{\lambda_i}(X) \right\}_{i=1}^{k}$ and $\left\{\phi_{\lambda_i}(X) \right\}_{i=1}^{k}$be sets of \textbf{k} unit eigenvectors.
		
		Further let $q$ be an observable function,
		\begin{equation}
			 q : M \to \mathbb{C}.
		\end{equation}
		
		If 
		\begin{equation}
			\displaystyle \min_{a} \Arrowvert \sum_{i=1}^{k} a_i \psi_{\lambda_i}(X) - q(X)  \Arrowvert \leq \min_{b} \Arrowvert \sum_{i=1}^{k} b_i \phi_{\lambda_i}(X) - q(X)  \Arrowvert
		\end{equation}
		
		then $\left\{\psi_{\lambda_i}(X) \right\}_{i=1}^{k}$ is a set of \textbf{k}-efficient Koopman eigenvectors
	\end{definition}

	\subsection{Computational Technique of Determining Koopman Eigen Pairs}
	
	Following definition of the eigenvector was rephrased from \cite{bollt2021geometric}.
	\begin{theorem}
		Let the Koopman eigenfunciton PDE is defined for an ODE $\dot{X} = F(X)$, with a flow $X(r) = \rho_r(X_0) : M \times \mathbb{R} \to M$. Assume a co-dimension-one \cite{abraham2008foundations} initial data manifold $\Lambda \subset M$ that is non-recurrent for some time epoch, $r\in \left[t_1, t_2\right]$ that contains $0$, and transverse to the flow, and let $\mathbb{U} = \cup_{t\in \left[t_1, t_2\right]} \rho_t(\Lambda)$ be the resulting non-recurrent closed domain. Furthermore, let $h : \Lambda \to \mathbb{C}$ be an initial data function then an Koopman eigen pair $(\lambda, \phi_{\lambda}(X))$ $\phi_{\lambda} : \mathbb{U} \to \mathbb{C}$ has the form,
		\begin{equation}
			\phi_{\lambda}(X) = h \circ s^*(X)e^{\lambda r*(X)}.
		\label{eq:eigen_vectors}
		\end{equation}
	\end{theorem}

	Now we present an computational explanation of the \textbf{Theorem 1}.  of \cite{bollt2021geometric}. This is the mathematical infra structure of the work presented in this paper.
	
	Suppose $\dot{X} = F(X)$ be an ODE satisfying the Koopman eigenvector PDE. Let the corresponding flow be given by $X(r) =  \rho_{r}(X_0) : M \times \mathbb{R} \to M$. Further let $\Lambda \subset M$ be a co-dimension-one data manifold. We should pick a time interval where the data values are non-recurrent. Let us call this $\Lambda$ a non-recurrent data manifold. Let us assume that the time span non-recurrent be $\left[t_1, t_2\right]$ and $0 \in \left[t_1, t_2\right]$.
	
	If $\Lambda$ is transverse to the flow, we can define $\mathbb{U}$ as follows. 
	
	\begin{equation}
		\mathbb{U} = \cup_{t\in \left[t_1, t_2\right]} \rho_t(\Lambda)
	\end{equation}
	
	That is $\mathbb{U}$ is the collection of all the trajectories starting from $\Lambda$. This includes reverse trajectories as well. Furthermore, let $h : \Lambda \to \mathbb{C}$ be an initial data function. Then an eigenpair $(\lambda, \phi_{\lambda}(X))$ satisfies the equation \ref{eq:eigen_vectors}. where,
	\begin{eqnarray}
		r^{*}(X) = \left\{t: \rho_{-t}(X)\cap \Lambda \neq \phi\right\}
	\end{eqnarray}
	and
	\begin{equation}
		s^{*} = s\circ \rho_{-r^{*}(X)}(X)
		\label{eq:s_star}
	\end{equation}

	Where $s$ is the parametrization of $\Lambda$. That is, $s$ is the data function available on $\Lambda$.
	
	The exercise of \cite{bollt2021geometric} was to introduce this formulation and introduce the algorithm to numerically compute $h$ and optimal $\lambda$ for a given stream of data.
	
	Following is a summary of the numerical algorithm described in \cite{bollt2021geometric}.
	
	Suppose $q \in L^2(\mathbb{U})$ be an arbitrary data function. $q : \mathbb{U} \to \mathbb{C}$. 
	Let,
	\begin{equation}
		\psi = argmin_{\lambda,h}\lVert\phi_{\lambda,h} - q\rVert_{L^2(\mathbb{U})}^2
		\label{eq:optimization}
	\end{equation}
	be the eigenvector that closely estimates $q$. $\psi : \mathbb{U} \to \mathbb{C}$. It is worth noting that the optimization is over both $\lambda$ and $h$. Further more,
	\begin{equation}
		argmax_{\lambda,h} \left(q, \phi_{\lambda,h}\right)= argmin_{\lambda,h}\lVert\phi_{\lambda,h} - q\rVert_{L^2(\mathbb{U})}^2
	\end{equation}
	where, 
	\begin{equation}
		\lVert f \rVert_{L^2(\mathbb{U})} = \int_{\mathbb{U}}f(X)dX
	\end{equation}
	and
	\begin{equation}
		(f,g) = \int_{\mathbb{U}}f(X)g(X)dX
	\end{equation}
	for any $f,g \in L^2(\mathbb{U})$.
	
	Let us explain an example scenario for a $2$-dimensional dynamical system. The algorithm can be generalized to larger dimensions limited only by the computational resources. Then $\Lambda$ becomes $1$-dimensional since it should be co-dimension $1$. let $s_0 <s_1<\dots <s_n$ be a uniform partition of $\Lambda$.
	The data function $h: \Lambda \to \mathbb{C}$ could be indexed accordingly as
	\begin{equation}
		h_i := h(s_i)
	\end{equation}

	We will consider a uniform grid of time as well. $r_0 = 0 < r_1 < \dots < r_m$ . That is $m+1$ number of equal time steps. Then the grid points can be identified as follows $s_{i,j} = \rho_{r_j} \circ x(s_i), 0\leq i \leq m, 0 \leq j \leq m$.
	
	Then in terms of the grid, the optimization problem in equation \ref{eq:optimization} for the function $\psi(x)$ represented on the grid $\psi \circ x(s_{i,j})$ is approximated by solving the finite rank least squares problem.
	
	\begin{equation}
		\tilde{\psi} = argmin_h \lVert \phi_{\lambda,h}\circ x(s_{i,j}) - q \circ x(s_{i,j})\rVert^2_{F} = argmin_{h_i} \sum_{j = 1}^{m} \sum_{i = 1}^{n}|e^{\lambda r_j} h _i - q_{i,j}|^2
	\end{equation}

 	On this grid both functions $q$ and $\tilde{psi}_{\lambda, h}$ can be sampled at $\rho_{r_j} \circ x(s_i)$.
 	
 	$q_{i,j} = q(s_{i,j})$ for $0\leq i \leq n, 0\leq m $. In \cite{bollt2021geometric}, for distance the Frobenius norm \cite{ackleh2009classical} is used. 
 	
 	We can rewrite the optimal initial data as an $n \times 1$ vector $h^{0}(\lambda) \in \mathbb{C}^n$. The problem restatement is as follows.
 	
 	Let us define the $m \times 1$ vector $E(\lambda)$ as $E(\lambda) = \left[e^{\lambda r_{0}} \quad e^{\lambda r_{2}} \dots e^{\lambda r_{m}}\right] \in \mathbb{C}^{m}$.
 	
 	Then $A(\lambda) = E(\lambda) \otimes I_n$ where $\otimes$ is the Kronecker product. Furthermore, let $b$ be the vector received by reshaping the matrix $q \in \mathbb{C}^{n \times m}$ into a vector of shape $mn \times 1$. That is $b = reshape(q, mn, 1)$.
 	
 	Then $h^0(\lambda)$ is computed by solving the least squares problem,
 	
 	\begin{equation}
 		h^0(\lambda) = argmin_h \lVert A(\lambda) h - b\rVert_{2}
 	\end{equation}
 
 	 It is noticed that $q = q(s_{i,j})$ is a grid of values. It follows that,
 	 
 	 \begin{equation}
 	 	p^{0}(\lambda) = A(\lambda) h^{0}(\lambda)
 	 \end{equation}
  
  Then $\tilde{\psi}^0_1 = reshape(p^0(\lambda), n,m)$.
  
  That is, in order to produce the optimal eigenvector, $p^0(\lambda)$ is being reshaped. 
  
  $\tilde{\psi}_1^0$ is the eigenvector computed by carrying out the computation once. This can be successively done to reduce the residuals.
  
  Therefore we can call the residue at the $k$th iteration $R_k$. 
  
  Then,
  \begin{equation}
  	R_k = b - \sum_{i = 1}^{k} p_{i}^0
  \end{equation}

	Then the successive eigen pairs can be computed by using,
	
	\begin{equation}
		(\lambda_{k+1}^0, \lambda_{k+1}^0) = argmin_{\lambda,h} \lVert A(\lambda)h - R_k \rVert_2
	\end{equation}

	\begin{equation}
		\tilde{\psi}^0_{k+1} = reshape(p_k^0,n,m)
	\end{equation}

 \section{Manifold Learning}\label{ref_Manifold Learning}

		In \cite{abraham2008foundations} and \cite{floryan2021charts} a manifold is described as a collection of charts which is called and atlas. Each chart maps a subset of the manifold to the \textit{Eucledean Space} with an invertible map. These subsests could have intersections which are not empty. The invertible maps should be consistent in the non empty intersections. Further more, a manifold is a topological space which is \textit{Locally Eucledean} \cite{abraham2008foundations}.

		\subsection{Manifold Hypothesis}
		The manifold hypothesis is the assumption that higher dimensional data lie on lower dimensional manifolds embedded within the higher dimensional space \cite{fefferman2016testing,narayanan2010sample}. This is due to the belief that high dimensional data could be generated from a low dimensional dynamical system.  Specifically in the context of a dynamical system, orbit data may be attracted to a stable invariant manifold, which occurs in certain systems where dissipation leads to the potential of a reduced order model.
   A collection of techniques developed to identify such a low dimensional manifold using the available high dimensional data is called manifold learning. For a comprehensive study of the manifold hypothesis and its mathematical consequences the reader is referred to \cite{fefferman2016testing}.

		\subsection{Manifold Learning}
		Among the handful of manifold learning techniques a few techniques such as isometric feature mapping \textit{ISOMAP}, Locally Linear Embeddings (\textit{LLE}), Laplacian Eigenmaps, Diffusion Maps, Nonlinear Principal Component Analysis stands out from the rest. While introducing manifold learning, \cite{izenman2012introduction} introducs it as an algorithmic technique for dimensionality reduction. In there paper \cite{tenenbaum2000global} introduces a geometric framework for nonlinear dimensionality reduction. The authors claim that their technique efficiently computes the global optimal and it is guaranteed to converge asymtotically to the true structure. For a comprehensive study of \textit{LLE} the reader is referred to \cite{roweis2000nonlinear}. For treatment of \textit{ISOMAP} algorithms and \textit{Laplacian Eigen Maps} the interested reader is referred to \cite{balasubramanian2002isomap} and \cite{belkin2001laplacian} respectively. 
		However this investigation is more interested about manifold learning using deep learning techniques. Following is a concise empirical study of deep learning based manifold learning in various engineering applications. For a more descriptive treatment on topological point of view of manifold learning refer to \cite{floryan2021charts}.
		
		\subsection{Mathematical Foundations of Manifold Learning}
		Dimensional reduction of the data is the main application of manifold learning. The data which is to be dimensionally reduced is assumed to be generated from a lower dimensional submanifold \cite{ma2020manifold}. A $d$-dimensional manifold $M$  has the property that for all $x \in M$ there is a neighborhood $X_x \subset M$ of $x$ such that $X_x$ is homeomorphic to an open set $\Theta_x \subset \mathbb{R}^d$ where $\mathbb{R}^d$ is the $d$-dimensional Euclidean space \cite{boothby2003introduction}.	Then $\Theta_x$ can be called a local coordinate of $X_x$ \cite{boothby2003introduction}. 
  
    Since $M = \bigcup\limits_{x \in M} X_x$ where $\{X_x | x \in M\}$ is an \textit{open covering} of $M$. We can find an homemorphism $\phi_x : X_x \rightarrow \Theta_x$ for all $x \in M$. Then for all $x \in M$, we have $\phi_x(X_x) = \Theta_x$. Therefore $\{\Theta_x | x \in M\}$ can be represented as $\{(X_x, \Theta_x) | x \in M \}$. The set of homeomorphisms $\{\phi_x\}_{x \in M}$ need not be unique \cite{boothby2003introduction}. This leads to the observation that learnt manifolds need not be unique. Suppose $X = \{X_1 , X_2, \dots , X_N \} \subset M \subset \mathbb{R}^{d'}$ be the set of observable data vectors. Let $Y = \{Y_1 , Y_2, \dots Y_N \}  \subset \mathbb{R}^d$ be the reduced dimensional vectors where $d' > d$. It can then be concluded that the observable data was generated by a $d$-dimensional dynamical system. Finding the homeomophism $\phi : X \rightarrow Y$ such that $\phi(X_i) = Y_i$ for $1 \leq i \leq N$ is the process of manifold learning.
		
		Even though the previously surveyed techniques can be used for manifold learning, this investigation is interested about deep learning techniques. The deep learning techniques are favourable in real time computing and the learnt deep neural network can be used to uncover the intrinsic local coordinates. This ability comes as a consequence of the results of this investigation.
		
		Following is a concise survey of applications of deep manifold learning. While introducing the phrase \textit{Deep Manifold Learning (DML)}  for deep learning for hidden less dimensional manifolds \cite{chen2017deep} represents a framework of DML for action recognition using convolutional Neural Networks. In their paper on image set classification \cite{wang2021symnet} introduces a technique for \textit{Riemannian} manifold deep learning. This paper also provides a summary of previous efforts on Riemmanian manifold learning. 		
		
	\subsection{Autoencoders}
\begin{figure}[h]
 		\centering
 		\includegraphics[width=0.8\linewidth]{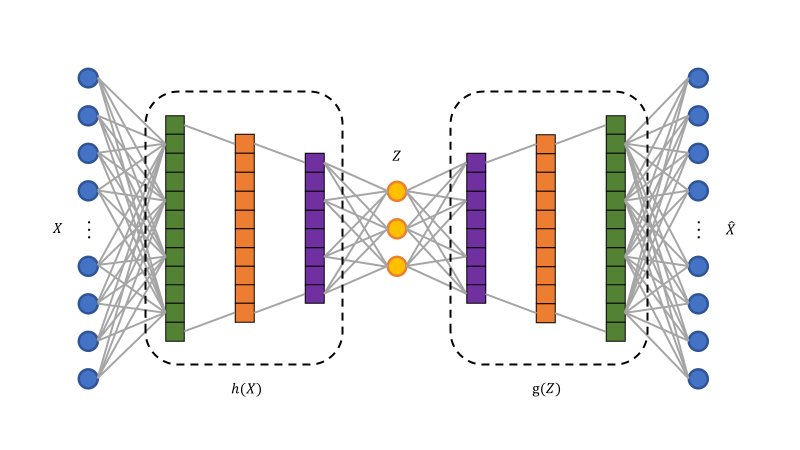}
 		\caption{Graphical depiction of a simplified deep autoencoder}
 		\label{fig:DAE}
 	\end{figure}

\begin{sloppypar}
Autoencoders have become a standout machine learning method to find a reduced order model in an invariant manifold within an artificial neural network framework, including in dynamical systems with stable reduced order models \cite{lee2020model, bakarji2022discovering}.  Figure \ref{fig:DAE} depicts a basic neural network topology of a deep autoencoder.
\end{sloppypar}
 
		Autoencoders (AE) were first proposed and formulated by the PhD thesis of Y. LeCun \cite{lecun1987phd}.These are traditionally used for dimensionality reduction \cite{zhai2018autoencoder}. Based on the formulation of \cite{gu2021autoencoder} the mathematical foundation of an AE can be viewed as a recursive formula \ref{eq:AE_recursive_formula}. The mathematical foundation of the neural networks can be found in \cite{bishop2006pattern}.

		\begin{equation}
			x^{(l)} = g\left(b^{(l-1)}+ W^{(l-1)}x^{(l-1)}\right)
			\label{eq:AE_recursive_formula}
		\end{equation}
	
		Let $K^{(l)}$ denote the number of neurons in each layer where, $l= 1,2,3,\dots,L$. Let the output of neuron $k$ in layer $l$ be $x_k^{(l)}$, and the vector of all outputs for this layer as $x^{(l)} = \left(x_{1}^{(l)},x_{2}^{(l)},x_{3}^{(l)},\dots,x_{K^{(l)}}^{(l)}\right)'$. In each node, a nonlinear activation function $g(\cdot)$ is applied further to the linear operation. In the recursive equation \ref{eq:AE_recursive_formula}, $W^{(l-1)}$ is a $K^{(l)} \times K^{(l-1)}$ matrix of weight parameters, and $b^{(l-1)}$ is a $K^{(l)} \times 1$ vector of bias parameters \cite{gu2021autoencoder}. For this work, "\textit{sigmoid}" and "\textit{ReLu}" \cite{karlik2011performance} were proven to be acceptable activation functions. We decided to do our experiments with \textit{sigmoid} consistently. 

    Autoencoders, in general, takes the form $Z = h(X)$ and $\widehat{X} = g(Z)$, where $\widehat{X}$ is the output of the neural network. Vector value $Z$ is known as the latent vector, which is the encoded/transformed form of $X$. Transformation $g$ suppose to invert the transformation $h$. Therefore, ideally, $g = h^{-1}$. The loss function of the neural network is a function of $\lVert \widehat{X} - X \rVert$. We used the mean squared error as the loss function,

    \begin{equation}
        \mathbb{L} = \frac{1}{N} \sum_{i = 1}^{N} \lVert \widehat{X} - X \rVert ^{2}
    \end{equation}

    Subsequently $\mathbb{L}$ is minimized using the gradient descent over a multitude of iterations.
	
	\subsection{Autoencoders for Reduced Order Modeling}
	
		Full order models require increased utilization of resources. The main goal of a reduced representation is to capture the essential physical features of a system and project them onto a lower-dimensional space or manifold in a way that preserves as much information as possible, while still allowing for meaningful comparisons to be made with a full-order model (FOM). This approach is therefore called reduced order modeling (ROM) \cite{pant2021deep}. There have been many attempts in developing ROMs for dynamical systems \cite{bai2002krylov, lucia2004reduced, kazantzis2010new}. Compiling a ROM for a complex dynamical system could be challenging task. Compiling a ROM for a complex dynamical system can be an exceedingly arduous task due to the inherent complexity of the system. The projection-based ROM is one of the earliest and most widely used techniques for reduced order modeling. This method involves transforming the original space into a lower-dimensional space using the governing PDEs (Partial Differential Equations) of the physical system \cite{san2018neural}.
		
		Proper Orthogonal Decomposition (POD) is a widely used projection-based ROM technique, particularly in the field of computational fluid dynamics. It involves decomposing the original data into a set of orthogonal modes that capture the dominant dynamics of the system. For those interested in delving deeper into the topic of POD, we recommend referring to the following seminal works: \cite{aubry1988dynamics, berkooz1993proper, holmes2012turbulence}. Principal Component Analysis (PCA) \cite{hotelling1933analysis} is a comparable method to POD which gathered momentum due to the availability of the computing resources. Empirical orthogonal functions \cite{lorenz1956empirical} Karhunen-Loeve expansion \cite{loeve1955probability} are also similar techniques for ROM. According to \cite{taira2017modal} these ROM techniques utilize the eigen-decomposition of the snapshot matrix using Singular Value Decomposition (SVD). 
		
		Dynamic Mode Decomposition (DMD) is a ROM technique which is popular in modeling physics based systems using there spatio-temporal coherent structures \cite{rowley2009spectral}. DMD is being used for a wide range of applications\cite{brunton2022data}.	Image-processing \cite{schmid2011applications}, Neuroscience \cite{brunton2016extracting}, and Robotics \cite{berger2014dynamic} are to name a few. Koopman operator theory (KOT) itself has direct connections to DMD. KOT was early used for characterization of the dynamics of the Hamiltonian functions \cite{koopman1931hamiltonian}. An Extended version of DMD (EDMD) was used in \cite{li2017extended} to decompose the KOT for the analysis of a non-linear dynamical system. A matching dynamical system is formed using the Koopman operator by decomposing the operator using EDMD in \cite{bollt2018matching}.

\section{Examples}\label{ref_Examples}

    Now we provide couple of worked examples of the developed theory. 

	\subsection{Example 1: Van der Pol Oscillator}
\begin{figure}[h]
 		\centering
 		\includegraphics[width=0.5\linewidth]{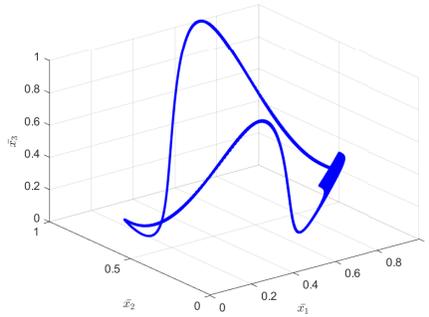}
 		\caption{Trajectories from the Van der Pol Oscillator projected on to the 3 dimensions and scaled}
 		\label{fig:VdP_trajectories_full_3d}
 	\end{figure}

\begin{figure}[h]
 		\centering
 		\includegraphics[width=0.5\linewidth]{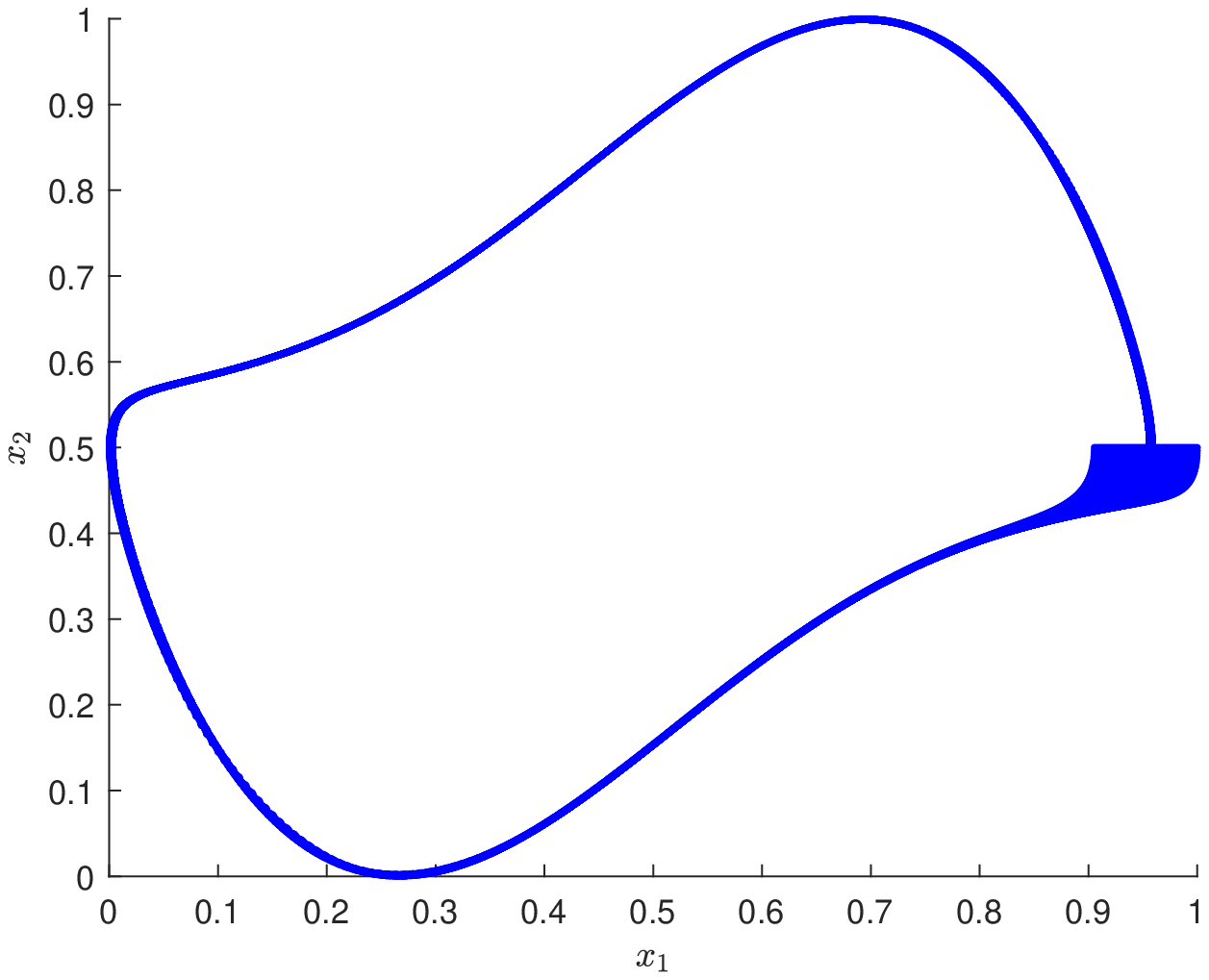}
 		\caption{Scaled trajectories from the Van der Pol Oscillator}
 		\label{fig:VdP_trajectories_full_2d}
 	\end{figure}

\begin{figure}[h]
 		\centering
 		\includegraphics[width=0.5\linewidth]{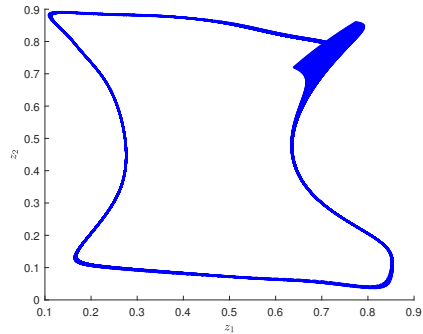}
 		\caption{Transformed trajectories from the Van der Pol Oscillator}
 		\label{fig:VdP_trajectories_full_transformed}
 	\end{figure}

	 We use the \textit{Van der Pol}(\textit{VdP}) oscillator as the first example for simulation. \textit{VdP} is an oscillator with non linear damping \cite{kanamaru2007van}. First step of the process is to simulate trajectories of the \textit{VdP}. In order to simulate geometric transformations using deep autoencoders, we considered a complete cycle of each trajectory. These trajectories are similar to what is depicted in figure \ref{fig:VdP_trajectories_full_2d}. The original trajectories lie on $2$ dimensional plane. Then we projected these trajectories to a $3$ dimensional space using the mapping,
	 \begin{equation}
	 	\begin{bmatrix} \bar{x}_1 \\ \bar{x}_2 \\ \bar{x}_3 \end{bmatrix}
	 	=
	 	\begin{bmatrix}
	 		x_1 \\
	 		x_2 \\
	 		x_1^2 + x_2^2
	 	\end{bmatrix}.
	 \end{equation}
 	This transforms $2d$ trajectories in to trajectories in $3d$. These projected trajectories are depicted in figure \ref{fig:VdP_trajectories_full_3d}. At this state the dynamical system is treated as a $3$ dimensional system.  The next step would be to reduce the order of the model to $2$. We did it using a deep-autoencoder(\textit{DAE}). Following is a brief description of the \textit{DAE} we tuned for this particular system. 
 	For implementing the deep neural network, we used \textit{Keras} which is based on \textit{Tensorflow}. Both of these were programmed using \textit{python} programing language. The \textit{DAE} had the latent vector size of $2$, consistent with our goal. We received optimal results when the first and second hidden layers had the sizes of $100$ nodes each. These two hidden layers are the encoder. The decoder, which is the neural network after the latent vector, has the same architecture as the encoder. Optimal results were received when \textit{Sigmoid} was used as the activation function. Comparably similar results were achieved when Rectified Linear Unit(\textit{Relu}) was used as the activation function. \textit{tanh} and \textit{softmax} functions did not comparably similar results. These two functions were rejected after visually inspecting the resulting reduced order trajectories. All the layers in the neural network were dense. As a normalized practice, the networks were trained for $2000$ epochs throughout the study. Before feeding into the neural network for training, each stream of data has to be scaled from $0$ to $1$. The scaled $2$ dimensional trajectories are depicted in figure \ref{fig:VdP_trajectories_full_2d}. There are many trajectories each with a different initial point. The nature of the \textit{VdP} is that they converge to a stable limit cycle. The range of initial points are depicted by the anomaly at around the points $(0.9,0.5)$ and $(1,0.5)$ in figure \ref{fig:VdP_trajectories_full_2d}. All the trajectories, initiated from the range of different initial points eventually converges to the stable limit cycle and stay on it. Therefore, it is sufficient to simulate one initial period. The transformed trajectory is depicted in figure \ref{fig:VdP_trajectories_full_transformed}. There is no guarantee that the transformed trajectory will always be this. Even with the original trajectories being the same the autoencoder might converge to a different local minimum of the error producing a different transformation. Finding a geometric transformation using an \textit{DAE} is a trial and error process. The number of hidden layers, number of nodes on each layer, number of epochs to train, the activation function may vary with the application at hand. Therefore it is, to be found and optimized over several iterations during the process of simulation.  
 	
 	\subsection{Reduced Order Modeling for Koopman Eigen Pair Computation}
\begin{figure}[h]
	 	\centering
	 	\includegraphics[width=0.5\linewidth]{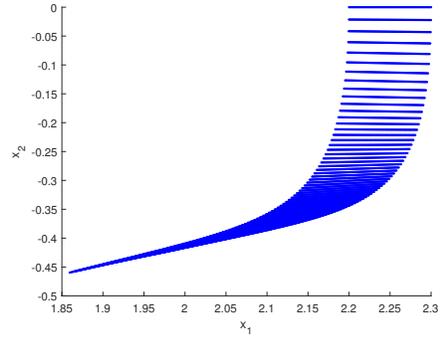}
	 	\caption{Trajectories from the Van der Pol Oscillator}
	 	\label{fig:VdP_trajectories}
	 \end{figure}

\begin{figure}[h]
	 	\centering
	 	\includegraphics[width=0.5\linewidth]{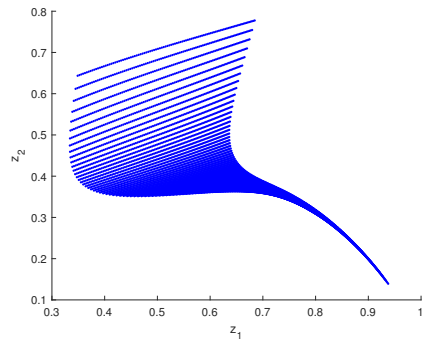}
	 	\caption{Trajectories in the transformed geometry}
	 	\label{fig:VdP_trajectories_encoded}
	 \end{figure}

 	Now we move into the analysis of the geometric transformation of the trajectory we used for \textit{Koopman} eigen pair computation. For the \textit{Koopman} analysis we do not need a full cycle of a trajectory. Therefore the trajectories are restricted to a shorter time period. The scaled trajectory is depicted in figure \ref{fig:VdP_trajectories}. Again this is projected to the $3$ dimensional space using the same transformation as in the complete trajectory case. The same neural network architecture as for the previous case was used to autoencode the set of shorter trajectories. The reduced order model (\textit{ROM}) trajectories are depicted in figure \ref{fig:VdP_trajectories_encoded}. A closer comparison between figures \ref{fig:VdP_trajectories_full_transformed} and \ref{fig:VdP_trajectories_encoded} will hint the interested reader that the two transformations are not equal to each other. For the \textit{DAE}, the same architecture and hyper-parameters were used as the full trajectory \textit{DAE}. This guarantees that accuracy is going to converge to a stable value. The next step is to use the \textit{ROM} to compute the \textit{Koopman} eigenvectors and eigenvalues.
 	
 	\subsection{Computing Koopman Eigenvectors with the \textit{ROM}}
\begin{figure}[h]
 		\centering
 		\includegraphics[width=0.8\linewidth]{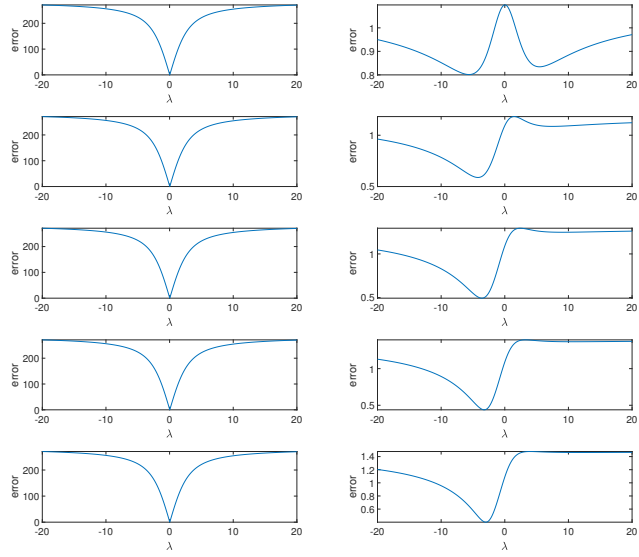}
 		\caption{Error vs $\lambda$ variation through the first $10$ modes}
 		\label{fig:VdP_error_vs_lambda}
 	\end{figure}

\begin{figure}[h]
 		\centering
 		\includegraphics[width=0.8\linewidth]{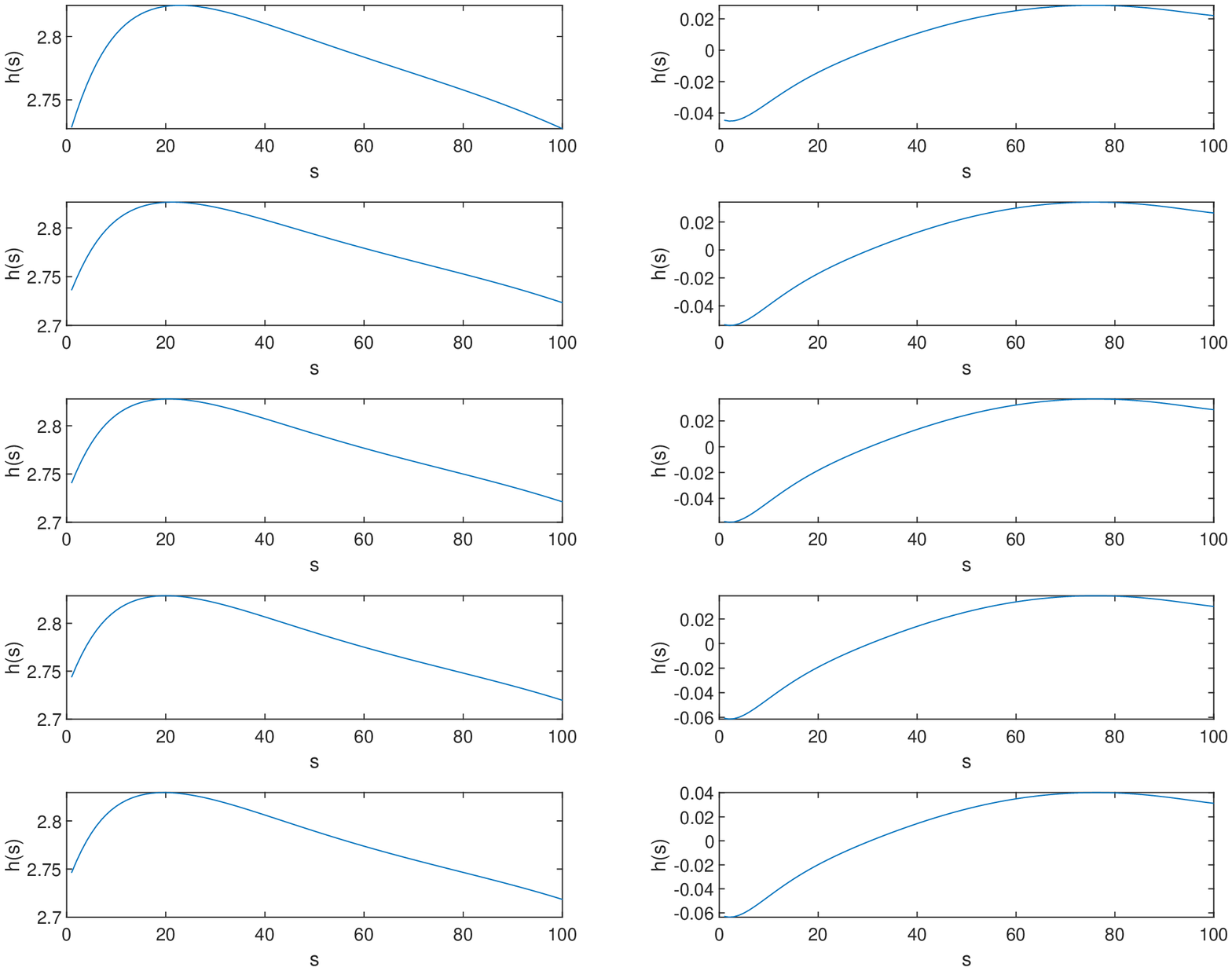}
 		\caption{$h(s)$ of the first $10$ modes}
 		\label{fig:VdP_h_vectors}
 	\end{figure}

\begin{figure}[h]
 		\centering
 		\includegraphics[width=0.6\linewidth]{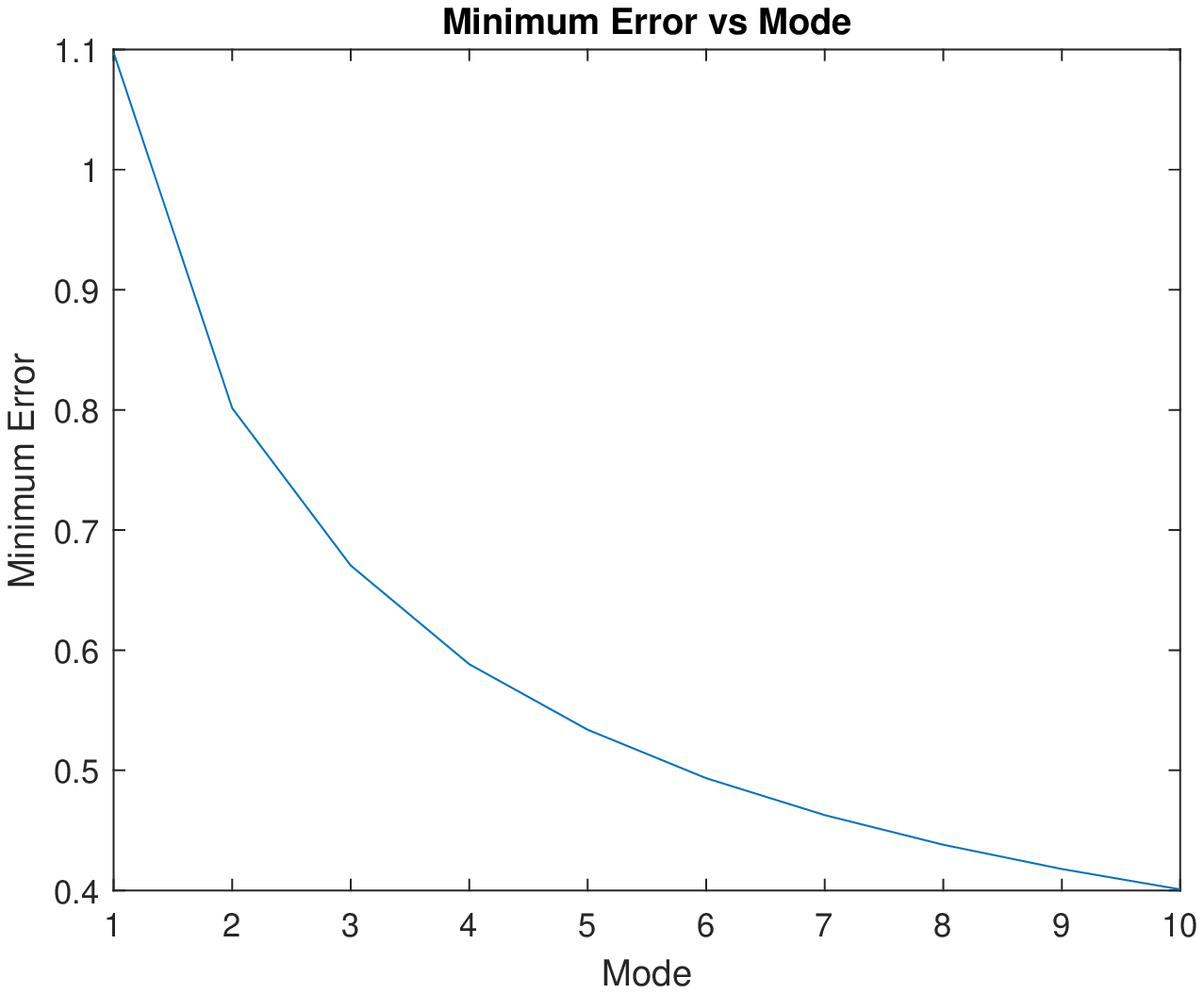}
 		\caption{Minimum Error vs. Mode}
 		\label{fig:min_error_vs_mode}
 	\end{figure}

 	This phase is to find \textit{Koopman} eigenvectors as given by the equation \ref{eq:eigen_vectors}. As given by the equation \ref{eq:s_star}, $s^*$ is derived as a consequence of the flow of the dynamical system. This leads to the fact that, we only need to determine the optimal eigenvalue $\lambda$ and the corresponding $h$ of the equation \ref{eq:eigen_vectors}. To be consistent with the investigation of \cite{bollt2021geometric} we chose the observation function to be $q = 3e^{-\frac{x_1^2 +x_2^2}{10}}$. When $x_i$s are replaced by $z_i$s we get,  $q = 3e^{-\frac{z_1^2 +z_2^2}{10}}$. In this case, instead of using $x_1$ and $x_2$ we use $z_1$ and $z_2$, the transformed coordinates. We computed the first $10$ modes or the first eigenvectors and their corresponding eigenvalues. For each mode the corresponding eigenvalue is the value of $\lambda$ which minimizes the error which is depicted in figure \ref{fig:VdP_error_vs_lambda}. The $h(s)$ vectors for the first $10$ modes are depicted in \ref{fig:VdP_h_vectors}. The absolute error between the sum of eigenvectors and the observation function decreases with the mode. This is depicted in figure \ref{fig:min_error_vs_mode}.

 	\begin{figure}[h]
 		\centering
 		\includegraphics[width=0.5\linewidth]{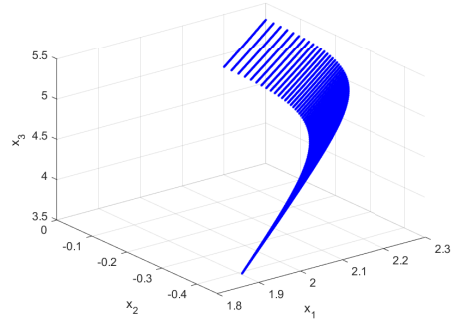}
 		\caption{Trajectories from the Van der Pol Oscillator projected to 3 dimensional space}
 		\label{fig:VdP_trajectories_3d}
 	\end{figure}
 
 	\subsubsection{eigenpairs with a Different Observation Function}
 	
\begin{figure}
 		\centering
 		\begin{subfigure}[b]{0.49\textwidth}
 			\centering
 			\includegraphics[width=\textwidth]{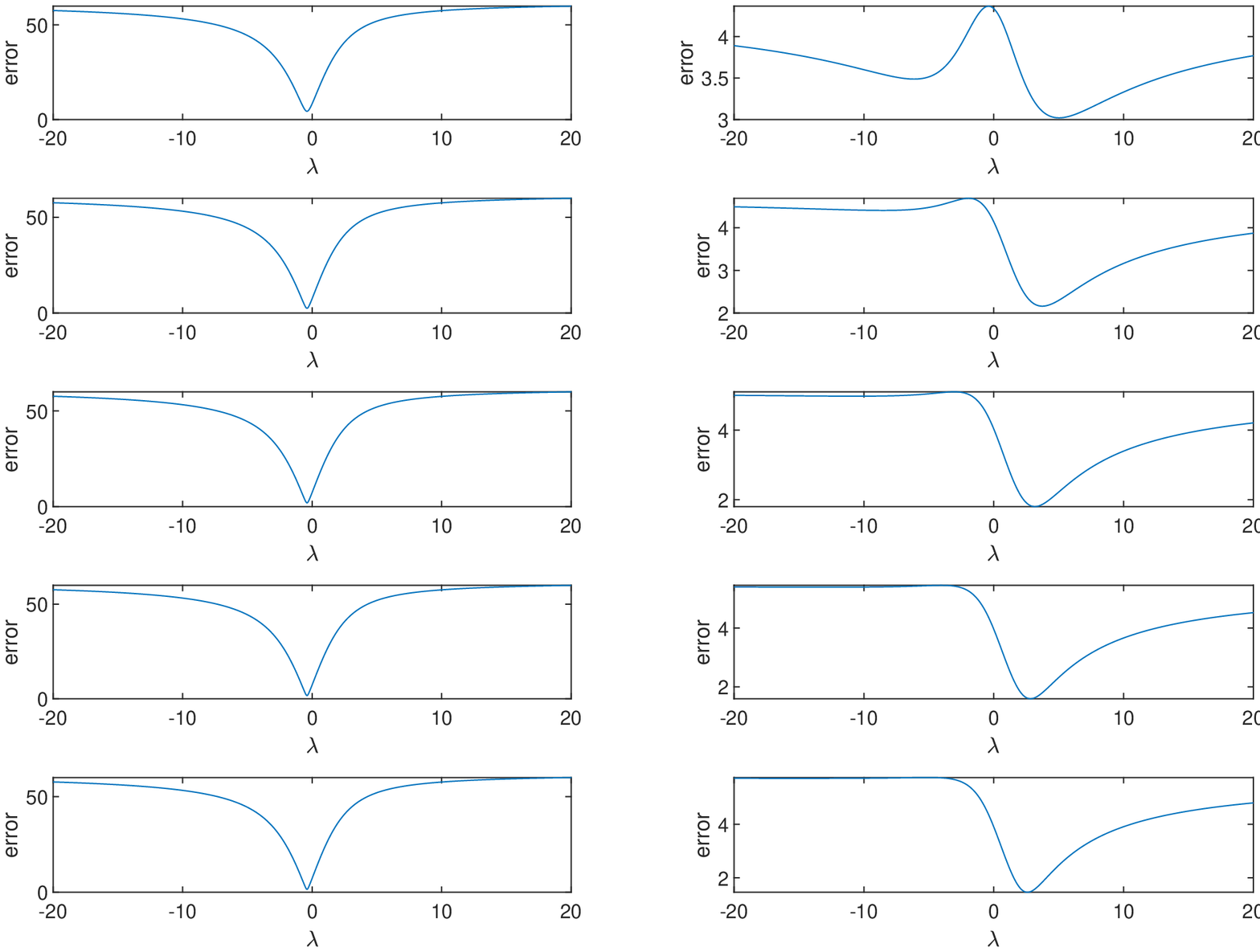}
 			\caption{Error vs $\lambda$ variation}
 			\label{fig:VdP_error_vs_lambda_q2}
 		\end{subfigure}
 		\hfill
 		\begin{subfigure}[b]{0.49\textwidth}
 			\centering
 			\includegraphics[width=\textwidth]{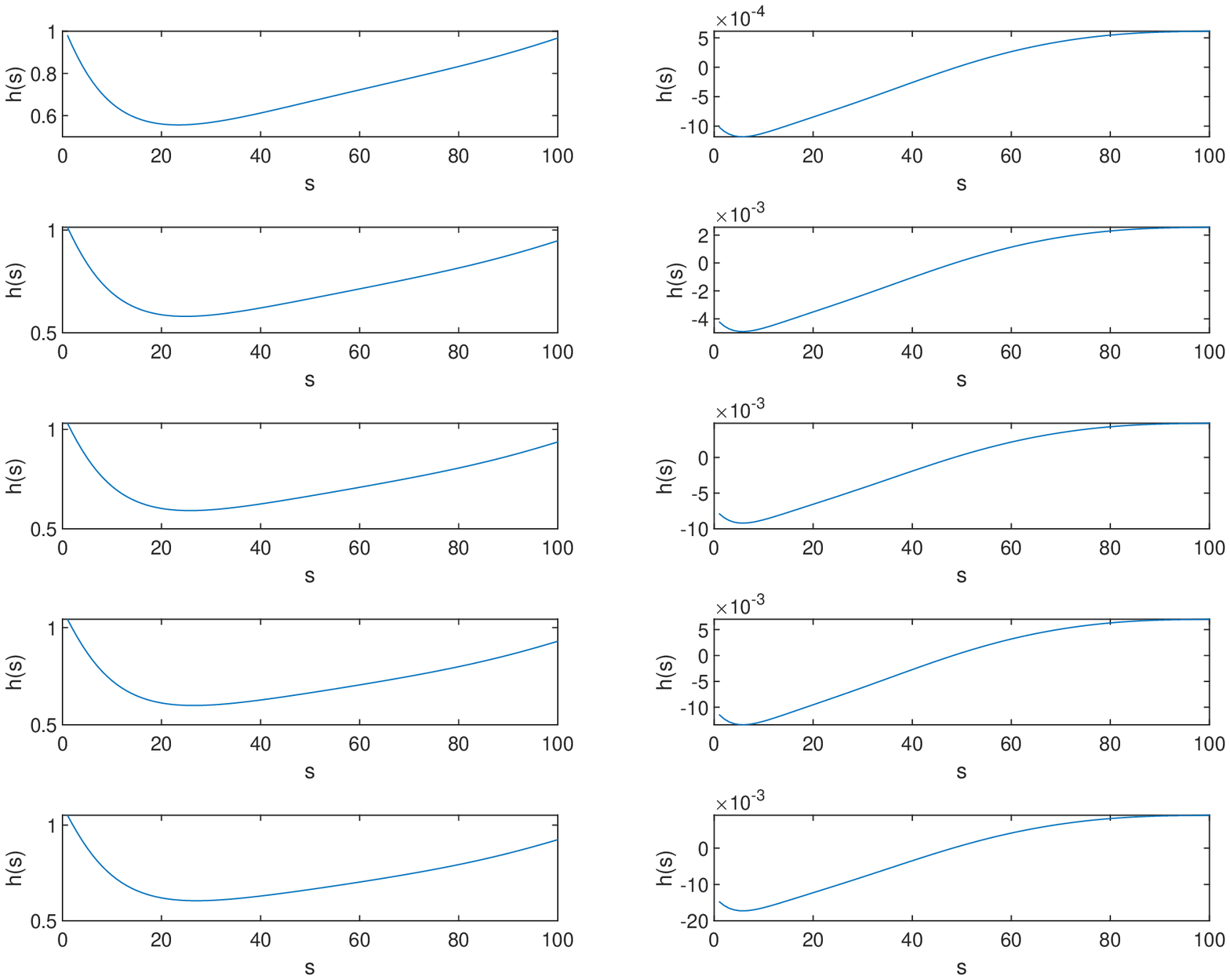}
 			\caption{$h(s)$}
 			\label{fig:VdP_h_vectors_q2}
 		\end{subfigure}
 		\hfill
 		
 		\caption{VdP trajectories with $q = z_1^2+z_2^2$}
 		\label{fig:VdP_lambadas_h_q2}
 	\end{figure}

\begin{figure}[h]
 	\centering
 	\includegraphics[width=0.6\linewidth]{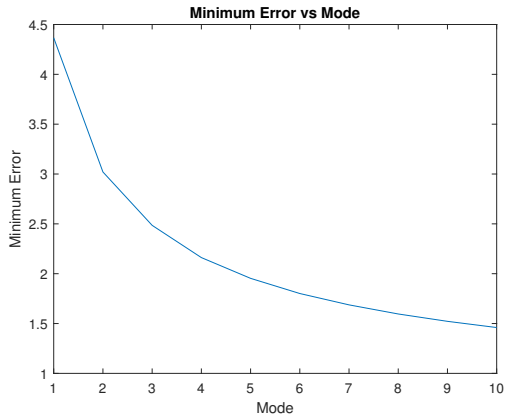}
 	\caption{Minimum Error vs. Mode for $q = z_1^2+z_2^2$}
 	\label{fig:min_error_vs_mode_q2}
 \end{figure}

 	Now we change the observation function to $q = z_1^2+z_2^2$. The variation of the error with $\lambda$ and the function $h(s)$ for the first $10$ modes are plotted in sub-figures, \ref{fig:VdP_error_vs_lambda_q2} and \ref{fig:VdP_h_vectors_q2} respectively. The corresponding minimum error vs mode graph is presented in figure \ref{fig:min_error_vs_mode_q2}
 
 \subsection{Example 2: Reaction Diffusion}
 	
 	Now we apply our algorithm to a $1$ dimensional reaction diffusion equation which appears in \cite{triandaf1997karhunen}.
 	
 	\begin{equation}
 		\begin{split}
	 		\frac{\partial u_1}{\partial t} &= D \frac{\partial^2 u_1}{\partial x^2} 	+ \frac{1}{\epsilon}\left(u_2 - f(u_1)\right)   \\
	 		\frac{\partial u_2}{\partial t} &= D \frac{\partial^2 u_2}{\partial x^2} 	-u_1 + \alpha
 		\end{split}
 		\label{eq:reaction_diffusion_RD}
 	\end{equation}
 
	 In the system \ref{eq:reaction_diffusion_RD}, $x\in \left[0,1\right]$ and it is subject to d\textit{Dirichilet} conditions,
 
	 \begin{equation}
	 	\begin{split}
	 		u_1\left(x=0, t=0\right)&= u_1\left(x=1, t=0\right) \\
	 		u_2\left(x=0, t=0\right)&= u_2\left(x=1, t=0\right).
	 	\end{split}
	 \label{eq:dirichilet_conditions}
	 \end{equation}
 
\begin{figure}[h]
 		\centering
 		\includegraphics[width=0.8\linewidth]{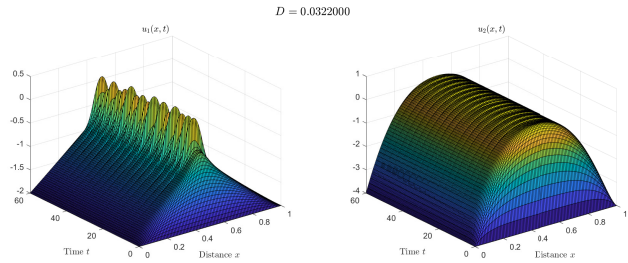}
 		\caption{Surface plots of $u_1$ and $u_2$ of Reaction-Diffusion PDE}
 		\label{fig:RD_u1_vs_u2}
 	\end{figure}

\begin{figure}[h]
 		\centering
 		\includegraphics[width=0.8\linewidth]{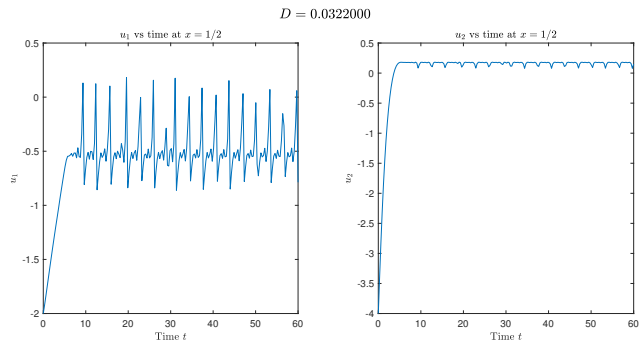}
 		\caption{$u_1$ and $u_2$ at $x=1/2$ of Reaction-Diffusion PDE}
 		\label{fig:RD_u1_vs_u2_mid_traj}
 	\end{figure}

 	We chose to fix $\epsilon = 0.01$ and $\alpha = 0.01$ for our experiments. The simulation we have presented, we used $D = 0.0322$. The surface plots for $u_1$ and $u_2$ for these parameters are presented in the figure \ref{fig:RD_u1_vs_u2} where $x \in \left[0,1\right]$ and $0 \leq t \leq 60$. It is noted that this reaction-diffusion system produces chaotic \cite{goldstein2002classical} behavior depending on the values of $\epsilon$ and $D$. The trajectories of $u_1$ and $u_2$ at $x= 1/2$ are graphed in figure \ref{fig:RD_u1_vs_u2_mid_traj}. Even though mathematically, there are unaccountably many trajectories between $x =0$ and $x = 1$, computationally this number is countably finite. We used $x$ to be \textit{linspace(0,1,100)}. This leads to the resulting simulation having $100$ trajectories. When there are $100$ trajectories each for $u_1$ and $u_2$ there are $200$ trajectories. Each trajectory from these becomes an input to the \textit{DAE}. In order the machine learn an reduced order model, the deep neural network had to have at least $4000$ nodes in the first hidden layer.  The network architecture with $(4000,4000)$ did not converge beyond $25\%$ accuracy. Increasing the number of nodes to $(10000,1000)$ did not improve the accuracy significantly. The available computing cloud resources were not sufficient to increase the number of nodes further. Therefore we had to look for other means of producing the reduced order model.

    The solution we found was to pre-process the data using other means. In particular, we used another embedding technique using time delayed snapshots of the trajectories.
 
 \subsection{Takens's Time Delay Embedding}

    Takens's theorem provides conditions where we can reproduce the dynamical system using sequential data from a single trajectory rather than processing all the trajectories in parallel. Since all the trajectories from the reaction-diffusion cannot be processed together in parallel with available computing resources, Takens's theorem can be used to reconstruct the dynamical system using time delay embedding with one trajectory at a time. 
 
 	\textit{Takens's} time delay embedding was introduced in the paper \cite{takens2006detecting}. This provides conditions under which a smooth attractor can be reconstructed from the observations made with an observable function. Suppose we have a $d$-dimensional dynamical system given by a state vector $x_t$ which is continues. Further, assume that we have one observable function $y(t)$ which is coupled to all components of $x_t$. Then a $k$ dimensional vector of observations can be created by considering $k$ time lagged observations with period $\tau$ of $y(t)$, that is $\left[\cdots,y_{t-2\tau},y_{t-\tau},y_{t},y_{t+\tau},y_{t+2\tau}, \cdots\right]$ so on and so forth. As $k \rightarrow \infty$ the system becomes deterministic and predictable. The \textit{Takens's} theorem states that the dynamics of the lagged vector becomes deterministic at a finite dimension. The finite dimension is given by $k < 2d+1$.
 	
 	Let us present the \textit{Takens's} embedding theorem more formally.
 	
 	\begin{theorem}[\textbf{Takens's Embedding Theorem}] Let $\dot{X} = f(X)$ be a dynamical system defined on the manifold $M$. $f : M \mapsto M$ and $f$ is smooth. Suppose that the dynamics $f$ has the strange attractor \cite{ruelle1971nature} $\mathbb{A}$ with Minkowski-Bouligand dimension \cite{falconer2004fractal} $d_{\mathbb{A}}$. Using the Whitney's embedding theorem \cite{adachi2012embeddings}, \cite{skopenkov2008embedding}, $\mathbb{A}$ can be embedded in $k$-dimensional Euclidean space with $k> 2d_{\mathbb{A}}$.
 	That is there is a diffeomorphism $\phi$ that maps $\mathbb{A}$ into $\mathbb{R}^k$ such that the derivative of $\phi$ has full rank. 
 	\label{thm:takens}		
 	\end{theorem}
 
 	Building on theorem \ref{thm:takens}, it is possible to construct a vector using only a single trajectory from the flow of trajectories of \ref{eq:reaction_diffusion_RD}.  In the work presented in this paper each of the trajectories $u_1$ and $u_2$ are time embedded into a vector of size $5$. It is further explained by using the trajectory at $x=1/2$. The values of the trajectory at $x=1/2$ are isolated for both $u_1$ and $u_2$. Since the values are discrete, $\tau$ is taken to be the time steps available. Let the time stamp we are considering to be $t=t_0$. Then the scalar $u_1\left(x = 1/2, t= t_0\right)$ is composed to be a vector given by, 
  $$\begin{bmatrix}
  u_1\left(x = 1/2, t_0 - 2\right)\\
  u_1\left(x = 1/2, t_0 - 1\right) \\  
  u_1\left(x = 1/2, t_0 \right)\\
  u_1\left(x = 1/2, t_0 + 1\right)\\
  u_1\left(x = 1/2, t_0 + 2\right)
  \end{bmatrix}.$$
  
  $\left[u_1\left(x = 1/2, t_0 - 2\right), u_1\left(x = 1/2, t_0 - 1\right), u_1\left(x = 1/2, t_0 \right),u_1\left(x = 1/2, t_0 + 1\right),u_1\left(x = 1/2, t_0 + 2\right)\right]^T$.
 	
 	When both $u_1$ and $u_2$ are considered, it becomes the following vector.
 		 	
			 	$$\begin{bmatrix}
			 		u_1\left(x = 1/2, t_0 - 2\right)\\ 
			 		u_1\left(x = 1/2, t_0 - 1\right)\\
			 		 u_1\left(x = 1/2, t_0 \right)\\
			 		 u_1\left(x = 1/2, t_0 + 1\right)\\
			 		 u_1\left(x = 1/2, t_0 + 2\right)\\
			 		 u_2\left(x = 1/2, t_0 - 2\right)\\ 
			 		 u_2\left(x = 1/2, t_0 - 1\right)\\
			 		 u_2\left(x = 1/2, t_0 \right)\\
			 		 u_2\left(x = 1/2, t_0 + 1\right)\\
			 		 u_2\left(x = 1/2, t_0 + 2\right)\\
			 	\end{bmatrix}$$

	Then the trajectories for different $x$ values are fed to the neural network serially.
	
	The first hidden layer had $3200$ nodes while the second hidden layer had $100$ nodes with a latent vector size of $6$ the network produced an accuracy $> 70\%$. Using this reduced order model the \textit{Koopman} eigenvectors were computed. 
	
	For the simulation of the eigen pair computations, $q = 3e^{-\frac{z_1^2 +z_2^2+z_3^2+z_4^2+z_5^2+z_6^2}{10}} $ was taken as the observations function. Similar to the previous experiment we have produced the graphs corresponding to the first $10$ modes. 
	
\begin{figure}[h]
 		\centering
 		\includegraphics[width=0.8\linewidth]{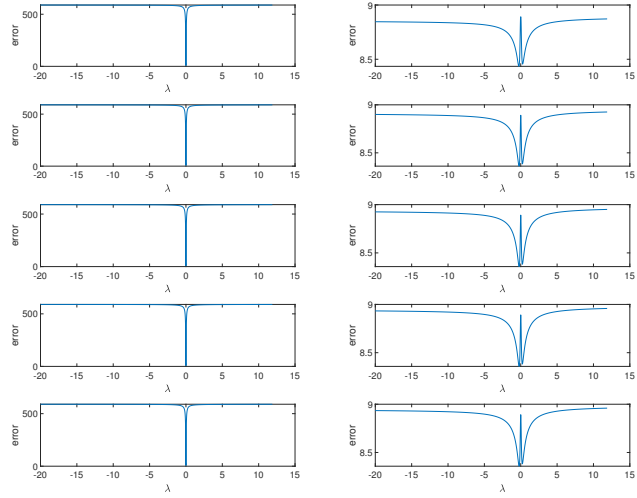}
 		\caption{Error vs $\lambda$ variation through the first $10$ modes Reaction-Diffusion PDE}
 		\label{fig:RD_error_vs_lambda}
 	\end{figure}

\begin{figure}[h]
 		\centering
 		\includegraphics[width=0.8\linewidth]{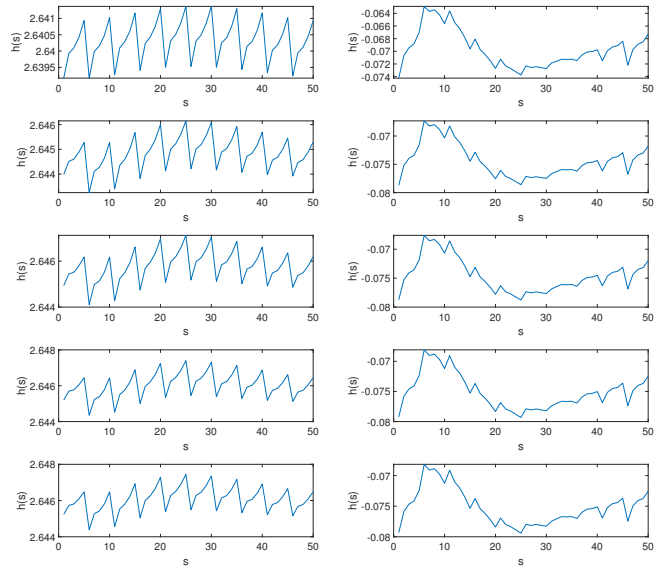}
 		\caption{$h(s)$ of the first $10$ modes Reaction-Diffusion PDE}
 		\label{fig:RD_h_vectors}
 	\end{figure}

\begin{figure}[h]
 		\centering
 		\includegraphics[width=0.6\linewidth]{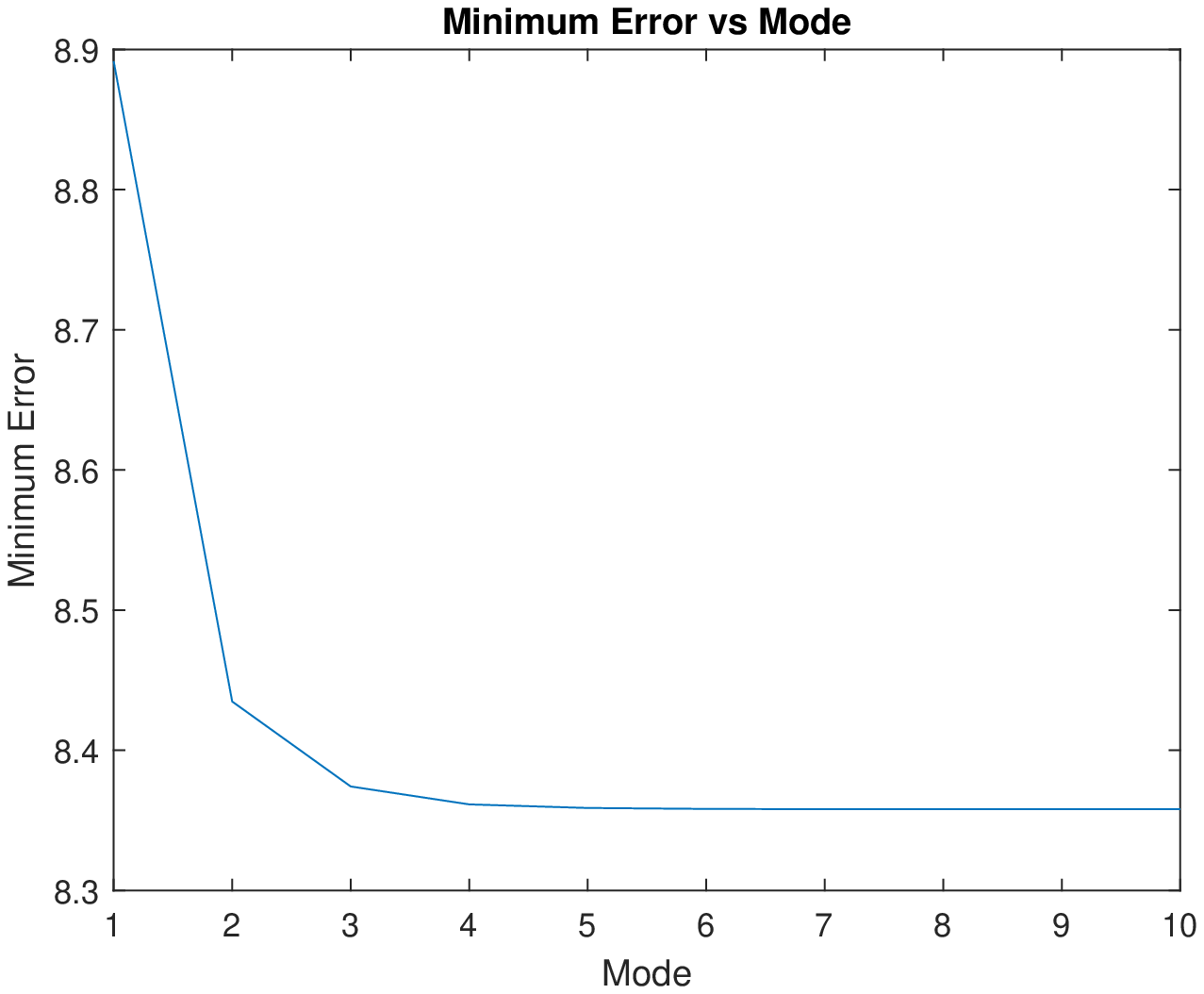}
 		\caption{Minimum Error vs. Mode Reaction-Diffusion PDE}
 		\label{fig:RD_min_error_vs_mode}
 	\end{figure}

	Figure \ref{fig:RD_error_vs_lambda} displays graphs that illustrate how the error changes with the eigenvalue of the initial $10$ modes. The binary nature of the graphs are observed again. While the graphs of the odd modes exhibit a similar shape, the graphs of the even modes display a shape that is different but yet similar to each other. Figure \ref{fig:RD_h_vectors} depicts the graphs of the computed $h(s)$ functions of the corresponding modes. It is visible that there is a "saw-tooth" nature to all the graphs even though there are two distinctive types of graphs. Figure \ref{fig:RD_min_error_vs_mode} depicts the error variation with mode. The absolute error is initially around $8.9$ but decreases to slightly below $8.4$, where it then remains constant.

\section{Conclusions and Further Work}\label{ref_Conclusions and Further Work}

In this paper we were trying to bring out deep autoencoding as a technique for uncovering simple geometries of complex dynamical systems. Autoencoder works as a geometric transformation. The transformed stream of data is then used to compute a dictionary of Koopman eigen pairs. The technique used to compute Koopman eigen pairs is different from the conventional technique and is introduced in the paper \cite{bollt2021geometric}. In that paper the eigen pairs are determined using the data streams as it is presented. In this paper, we improve on that and compute the eigen pairs using the transformed geometry.

The transformation minimizes the vector size of the input data by compressing the data driven dynamical system into a low dimensional manifold. This reduces the required computational resources to compute the eigen pairs. The transformed data stream is then a reduced order model of the full order dynamical system.

When this \textit{ROM} is formed using an autoencoder there is no technique to make this model according to our need. This is because the autoencoder is an unsupervised learning technique. To overcome this, a non-trivial loss function can be introduced. This loss function could partly depend on the latent variables and partly on the output of the neural network. Even though this is still unsupervised learning, the \textit{ROM} can be governed using a non trivial loss function. This will, in turn, lead to less error in Koopman eigen pair calculations.

\bibliography{NRef}
\bibliographystyle{amsplain}

\end{document}